\documentstyle[12pt, page1]{article}
\newcommand{\sect}[1]{\section{#1}\setcounter{equation}{0}}

\font\mbn=msbm10 scaled \magstep1
\font\mbs=msbm7 scaled \magstep1
\font\mbss=msbm5 scaled \magstep1
\newfam\mbff
\textfont\mbff=\mbn
\scriptfont\mbff=\mbs
\scriptscriptfont\mbff=\mbss
\def\mbf{\fam\mbff}
\def\Re{{\mbf R}}

\def\Z{{\mbf Z}}
\def\Co{{\mbf C}}
\def\To{{\mbf T}}
\def\P{{\mbf P}}
\def\Di{{\mbf D}}
\def\B{{\mbf B}}
\newtheorem{Th}{Theorem}[section]
\newtheorem{Lm}[Th]{Lemma}
\newtheorem{C}[Th]{Corollary}

\newtheorem{E}[Th]{Example}
\author{Alexander Brudnyi\thanks{Research supported in part by NSERC.
}\\
Department of Mathematics and Statistics\\
University of Calgary\\
Calgary, Canada}
\title{A Uniqueness Property for $H^{\infty}$ on Coverings
of Projective Manifolds}
\date{}
\begin{document}
\maketitle
\begin{abstract}
{Let $Y$ be a regular covering of a complex projective manifold
$M\hookrightarrow\Co\P^{N}$ of dimension $n\geq 2$. Let $C$ be intersection
with $M$ of at most $n-1$ generic hypersurfaces of degree $d$ in 
$\Co\P^{N}$. The preimage $X$ of $C$ in $Y$ is a connected submanifold.
Let $H^{\infty}(Y)$ and $H^{\infty}(X)$ be the Banach spaces of bounded
holomorphic functions on $Y$ and $X$ in the corresponding 
supremum norms. We prove that the restriction 
$H^{\infty}(Y)\longrightarrow H^{\infty}(X)$ is an isometry for $d$ large
enough. This answers the question posed in [L].\\}
\end{abstract}
\sect{\hspace*{-1em}. Introduction and Formulation of the Result.}
{\bf 1.1. An Extension Theorem.} Let $M$ be a complex projective manifold of
\penalty-10000
dimension $n\geq 2$ with a K\"{a}hler form $\omega$ and let $L$ be a positive
line bundle on $M$ with canonical connection $\nabla$ and curvature 
$\Theta$ in a hermitian metric $h$. Let $C$ be the common zero locus of
holomorphic sections $s_{1},...,s_{k}$, $k<n$, of $L$ over $M$ which, in a
trivialization, can be completed to a set of local coordinates at each point
$C$. Then $C$ is a (possibly disconnected) $k$-dimensional submanifold of 
$M$ which will be referred to as an {\em $L$-submanifold of $M$}. 
Let $\pi:Y_{G}\longrightarrow M$ be a regular covering of $M$ with a 
transformation group $G$ and $X_{G}=\pi^{-1}(C)$.  We denote
the pullbacks to $Y_{G}$ of $\omega$ and $\Theta$ by the same letters.
\begin{E}\label{e0}
{\rm If $L$ is very ample, then it is pullback of
the hyperplane bundle by an embedding of $M$ into some projective space
$\Co\P^{N}$. Further, zero loci of holomorphic sections of $L$ are
hyperplane sections of $M$. By Bertini's theorem, the generic linear
subspace of codimension $n-k$, $k<n$, intersects $M$ transversely in
a smooth manifold $C$ of dimension $k$, and by the Lefschetz hyperplane 
theorem, $C$ is connected and the induced homomorphism 
$\pi_{1}(C)\longrightarrow \pi_{1}(M)$ of fundamental groups is surjective.
Thus in this case $X_{G}\subset Y_{G}$ is a connected submanifold.}
\end{E}
Further, let $dist(\cdot,\cdot)$ be the distance on $Y_{G}$ induced by 
$\omega$.
Consider a function $\phi: Y_{G}\longrightarrow\Re$ such that $d\phi$ is 
bounded, i.e. 
$$
|\phi(x)-\phi(y)|\leq a\cdot dist(x,y)\  \ \ \ 
{\rm for\ some}\ \ \ a>0.
$$
By ${\cal O}_{\phi}(X_{G})$ we denote the vector space of 
holomorphic functions on $X_{G}$ such that $|f|^{2}e^{-\phi}$ is integrable 
on $X_{G}$ with respect to the 
volume form of the induced K\"{a}hler metric on $X_{G}$. This is a Hilbert
space with respect to the inner product
$$
(f,g)\mapsto\int_{X_{G}}f\overline{g}e^{-\phi}\omega^{k}\ .
$$
We define ${\cal O}_{\phi}(Y_{G})$ similarly. By $|\cdot|_{\phi,X_{G}}$
and $|\cdot|_{\phi,Y_{G}}$ we denote the corresponding norms. It
was shown in [L] that the restriction determines a continuous linear map
$$
\rho: {\cal O}_{\phi}(Y_{G})\longrightarrow {\cal O}_{\phi}(X_{G}),\ \ \
f\mapsto f|_{X_{G}}\ .
$$
The following remarkable result was proved by L\'{a}russon [L, Th.1.2].
\begin{Th}\label{la}
Suppose
$$
\Theta\geq i\partial\overline{\partial}\phi+\epsilon\omega
$$
for some $\epsilon>0$ in the sense of Nakano. Then $\rho$ is an 
isomorphism.
\end{Th}
{\bf 1.2. Formulation of the Main Result.}
An important example of a function $\phi$ as above is obtained by smoothing
the distance $\delta$ from a fixed point $o$ in $Y_{G}$. By a result of Napier
[N], there is a smooth function $\tau$ on $Y_{G}$ such that
\begin{description}
\item[{\rm (A)}] $c_{1}\delta\leq\tau\leq c_{2}\delta+c_{3}$ for some
$c_{1},c_{2},c_{3}>0$,
\item[{\rm (B)}] $d\tau$ is bounded, and
\item[{\rm (C)}] $i\partial\overline\partial\tau$ is bounded.
\end{description}
Furthermore, by (A) and since the curvature of $Y_{G}$ is bounded below,
there is $c>0$ such that $e^{-c\tau}$ is integrable on $Y_{G}$. Then
$e^{-c\tau}$ is also integrable on $X_{G}$. We set 
\begin{equation}\label{A}
A:=\frac{cc_{2}}{c_{1}}\ .
\end{equation}

Let $\tilde L$ be any positive line bundle on $M$ with curvature 
$\tilde\Theta$. By
(C) there is a non-negative integer $m_{0}$ such that for any integer 
$m>m_{0}$
\begin{equation}\label{e1}
m\tilde\Theta> i\partial\overline\partial(A\tau)\ .
\end{equation}
We set $L:=\tilde L^{\otimes m}$. 
Then L\'{a}russon's theorem holds for 
coverings $X_{G}:=\pi^{-1}(C)\subset Y_{G}$ of $L$-submanifolds $C\subset M$
with $\phi:=A\tau$ and with $\phi:=c\tau$ (because $A\geq c$).

Let $H^{\infty}(Y_{G}),\ H^{\infty}(X_{G})$ be the Banach 
spaces of bounded holomorphic functions on $Y_{G}$ and $X_{G}$ in
the corresponding supremum norms.
\begin{Th}\label{te1} 
The map $\rho: H^{\infty}(Y_{G})\longrightarrow H^{\infty}(X_{G})$, 
$f\mapsto f|_{X_{G}}$, is an isometry.
\end{Th}
{\bf 1.3. Corollaries and Examples.} Let $X$ be a complex manifold and
$H^{\infty}(X)$ be the Banach algebra
(in the supremum norm) of bounded holomorphic functions on $X$. The maximal
ideal space ${\cal M}={\cal M}(H^{\infty}(X))$ is the set of all nontrivial 
linear multiplicative functionals on $H^{\infty}(X)$. The norm of any 
$\phi\in {\cal M}$ is $\leq 1$ and so
${\cal M}$ is embedded into the unit ball of the dual space 
$(H^{\infty}(X))^{*}$. Then
${\cal M}$ is a compact Hausdorff space in the weak $*$ topology 
induced by $(H^{\infty}(X))^{*}$ (i.e. the {\em Gelfand topology}).
Further, there is a continuous map $i:X\longrightarrow {\cal M}$ 
taking $x\in X$ to the evaluation homomorphism $f\mapsto f(x)$. This map is
an emebedding if $H^{\infty}(X)$ separates points of $X$.
Recall also that the complement to the closure of $i(X)$ in 
${\cal M}$ is called {\em the corona}. The {\em corona problem} is to 
determine those $X$ for which the corona is empty.  
For example, according to Carleson's
celebrated Corona Theorem [C] this is true for $X$ being
the open unit disk $\Di\subset\Co$. Also there are non-planar Riemann 
surfaces for
which the corona is non-trivial (see e.g. [JM], [G], [BD] and references 
therein).
The general problem for planar domains is still
open, as is the problem in several variables for the ball and polydisk.
In [L, Th. 2.1] L\'{a}russon discovered simplest examples of Riemann 
surfaces with big corona. Namely, using his Theorem \ref{la} he proved that
if $Y_{G}\subset\Co^{n}$ is a bounded domain and $X_{G}\subset Y_{G}$ is a
Riemann surface satisfying assumptions of Theorem \ref{te1} then the natural 
map $i:X_{G}\hookrightarrow {\cal M}(H^{\infty}(X_{G}))$ extends to an 
embedding $Y_{G}\hookrightarrow  {\cal M}(H^{\infty}(X_{G}))$. The 
next statement is an extension of his result.
\begin{C}\label{cor1}
Under the assumptions of Theorem \ref{te1}, the transpose map \penalty-10000
$\rho^{*}:{\cal M}(H^{\infty}(X_{G}))\longrightarrow 
{\cal M}(H^{\infty}(Y_{G}))$,\ $\phi\mapsto\phi\circ\rho$, is a homeomorphism.
\end{C}
This follows from the fact that $\rho: H^{\infty}(Y_{G})\longrightarrow
H^{\infty}(X_{G})$ is an isometry of Banach algebras.\ \ \ \ \ $\Box$
\begin{E}\label{e1}
{\rm (1) (The references for this example are in [L, Sect. 4].)
Let $M$ be a projective manifold covered by the unit ball 
$\B\subset\Co^{n}$ with a positive line bundle $L$ with curvature $\Theta$,
and $X\subset\B$ be the preimage of an $L$-submanifold $C\subset M$. Let
$\delta$ be the distance from the origin in the Bergman metric of $\B$.
By $\omega$ we denote the K\"{a}hler form of the Bergman metric.
It was shown in [L, Sect. 4] that there is a nonnegative 
function $\tau$ on $\B$ such that 
$i\partial\overline{\partial}\tau=\omega$, $d\tau$ is bounded, and
$$
\sqrt{n+1}\delta\leq\tau\leq\sqrt{n+1}\delta+(n+1)\log\ 2 .
$$
Moreover, 
$$
\int_{\B}e^{-c\tau}\omega^{n}<\infty\ \ \ \ {\rm if\ and\ only\ if}\ \ \
c>\frac{2n}{n+1}\ .
$$
Applying Theorem \ref{te1} (with $c_{2}=c_{1}=\sqrt{n+1}$) we obtain that
{\em $\rho:H^{\infty}(\B)\rightarrow H^{\infty}(X)$ is an isometry
if $\Theta>\frac{2n}{n+1}\omega$.} This holds for instance if 
$L=K^{\otimes m}$ with $m\geq 2$ where $K$ is the canonical bundle of $M$.\\
(2) Let $S$ be a compact complex curve of genus 
$g\geq 2$ and $\Co\To$ be a one-dimensional complex torus. 
Consider an $L$-curve $C$ in 
$M:=S\times\Co\To$ with a very 
ample $L$ satisfying the assumptions of Theorem \ref{te1}. Let 
$\pi:\Di\times\Co\longrightarrow M$ be the universal covering. Then 
Theorem \ref{te1} is valid for the
connected curve $X:=\pi^{-1}(C)\subset\Di\times\Co$. This implies that 
any $f\in H^{\infty}(X)$ is constant on each 
$S_{y}:=(\{y\}\times\Co)\cap X$, $y\in\Di$. Note that $S_{y}$ is union of 
the orbits of some $z_{iy}\in X$, $i=1,...,k$, under the natural action of 
the group $\pi_{1}(\Co\To)\ (\cong\Z\oplus\Z)$ on $\Di\times\Co$.}
\end{E}
\sect{\hspace*{-1em}. Proof of Theorem 1.3.}
Let us fix a fundamental compact $K$ of the action of $G$ on 
$Y_{G}$, i.e., $Y_{G}=\cup_{g\in G}\ g(K)$. Consider finite covers
${\cal U}=(U_{i})$ and ${\cal V}=(V_{j})$ of $K$ by 
compact coordinate polydisks such that each $V_{j}$ belongs to the interior 
of some $U_{i_{j}}$.
\begin{Lm}\label{le1} 
Let $f$ be a holomorphic function defined in an open
neighbourhood $O$ of $\cup_{i}\ U_{i}$. Assume that
$$
\int_{O}|f|^{2}\omega^{n}=B<\infty\ .
$$
Then there is a constant $b>0$ (depending on ${\cal U}$ and ${\cal V}$
only) such that
\begin{equation}\label{e2}
\max_{K}|f|\leq b\sqrt{B}\ .
\end{equation}
\end{Lm}
The proof of the lemma is the consequence of the following facts:\\
(a) after the identification of $U_{i_{j}}$ with the closed unit polydisk 
$D$ and of 
$V_{j}$ with a compact subset $D_{j}\subset D$,
the volume form $\omega^{n}$ restricted to each $U_{i_{j}}$ is equivalent
to the Euclidean volume form
$do:=dz_{1}\wedge d\overline{z}_{1}\wedge\dots\wedge
dz_{n}\wedge d\overline{z}_{n}$ ;\\
(b) the Bergman inequality (see [GR, Ch.6, Th.1.3])
$$
\max_{D_{j}}|f|\leq\frac{(\sqrt{n})^{n}}{(\sqrt{\pi}d)^{n}}\cdot
\left(\int_{D}|f|^{2}do\right)^{1/2},
$$
where $d$ is the Euclidean distance from $D_{j}$ to the boundary of $D$.\\
We leave the details to the reader.\ \ \ \ \  $\Box$

Let $dist(\cdot,\cdot)$ be the distance on $Y_{G}$ in the metric induced
by $\omega$. In particular, $\delta(x):=dist(x,o)$. Since $\omega$ is 
invariant with respect to the action of $G$ we also have 
$dist(g(x),g(y))=dist(x,y)$ for any $g\in G$. From 
inequalities (A) for $\tau$ and the triangle inequality
for the distance we obtain
\begin{equation}\label{e3}
\begin{array}{c}
\tau(g(x))\geq c_{1}dist(g(x),o)\geq c_{1}[dist(g(x),g(o))-dist(g(o),o)]=\\
\\
c_{1}[dist(x,o)-dist(g(o),o)]\geq
(c_{1}/c_{2})\tau(x)-(c_{1}c_{3}/c_{2})-
c_{1}\delta(g(o))\ .
\end{array}
\end{equation}
Further, if $x\in K$ then 
\begin{equation}\label{e4}
a_{1}\leq\tau(x)\leq a_{2}\ \  {\rm for\ some}\ \ a_{1},a_{2}>0\ .
\end{equation}

Below by $|\cdot|_{\infty, X_{G}}$ and
$|\cdot|_{\infty, Y_{G}}$ we denote the corresponding $H^{\infty}$-norms.
Let
$f\in H^{\infty}(X_{G})$. 
Then $f\in {\cal O}_{A\tau}(X_{G})\cap {\cal O}_{c\tau}(X_{G})$ and 
there is $a_{3}>0$ such that
$$
\max\{|f|_{A\tau, X_{G}}, |f|_{c\tau, X_{G}}\}\leq  a_{3}\sup_{X_{G}}|f|:=
a_{3}|f|_{\infty, X_{G}}\ .
$$
By Theorem \ref{la}, there is a unique
$\tilde f\in {\cal O}_{A\tau}(Y_{G})\cap {\cal O}_{c\tau}(Y_{G})$ such that
$\tilde f|_{X_{G}}=f$ and
$$
\max\{|\tilde f|_{A\tau, Y_{G}}, |\tilde f|_{c\tau, Y_{G}}\}\leq a_{4}
\max\{|f|_{A\tau, X_{G}}, |f|_{c\tau, X_{G}}\}\ \ {\rm for\ some}\
\ a_{4}>0 .
$$
Combining these inequalities with (\ref{e4}) and (\ref{e2}) yields
$$
\max_{K}|\tilde f|\leq a_{5}|f|_{\infty, X_{G}},
$$
with some $a_{5}>0$ depending on $X_{G},Y_{G}$ only. For a fixed $g\in G$
consider $(g^{*}f)(z):=f(g(z))$. As above, there is a unique function
$\tilde f_{g}\in {\cal O}_{A\tau}(Y_{G})\cap {\cal O}_{c\tau}(Y_{G})$,
$\tilde f_{g}|_{X_{G}}=g^{*}f$,
such that
$$
\max_{K}|\tilde f_{g}|\leq a_{5}|f|_{\infty, X_{G}},
$$
But according to (\ref{e3}) and (\ref{A}), function
$(g^{*}\tilde f)(z):=\tilde f(g(z))$ belongs to ${\cal O}_{A\tau}(Y_{G})$
and $(g^{*}\tilde f-\tilde f_{g})|_{X_{G}}\equiv 0$. Thus by 
Theorem \ref{la} we have $\tilde f_{g}=g^{*}\tilde f$. Since $K$ is the
fundamental compact, the above inequality for each $\tilde f_{g}$ implies that
\begin{equation}\label{e5}
|\tilde f|_{\infty, Y_{G}}\leq a_{5}|f|_{\infty, X_{G}}\ .
\end{equation}
We will prove now that $a_{5}=1$ which gives the required statement.
Indeed, the same arguments as above show that for any integer $n\geq1$ the
function $(\tilde f)^{n}$ is the unique extension of $f^{n}$ satisfying 
(\ref{e5}):
$$
|(\tilde f)^{n}|_{\infty, Y_{G}}\leq a_{5}|f^{n}|_{\infty, X_{G}}\ .
$$
Thus 
$$
|\tilde f|_{\infty, Y_{G}}\leq 
\lim_{n\to\infty}(a_{5})^{1/n}|f|_{\infty, X_{G}}=|f|_{\infty, X_{G}}
$$
The proof of the theorem is complete.\ \ \ \ \ \ $\Box$


\begin{thebibliography}{   }
\bibitem[BD]{BD}
D. E. Barett and J. Diller, A new construction of Riemann surfaces with
corona. J. Geom. Anal. {\bf 8} (1998), 341-347.
\bibitem[C]{C}
L. Carleson, Interpolation of bounded analytic functions and the corona
problem. Ann. of Math. {\bf 76} (1962), 547-559.
\bibitem[G]{G}
T. W. Gamelin, Uniform algebras and Jensen measures. London Math. Soc.
Lecture Notes Series {\bf 32}, Cambridge University Press, 1978.
\bibitem[GR]{GR}
H. Grauert and R. Remmert, Theorie der Steinschen R\"{a}ume. 
Springer-Verlag, Berlin, 1977.
\bibitem[JM]{JM}
P. Jones and D. Marshall, Critical points of Green's functions, harmonic 
measure and the corona theorem. Ark. Mat. {\bf 23} no.2 (1985), 281-314.
\bibitem[L]{L}
F. L\'{a}russon, Holomorphic functions of slow growth on nested covering
spaces of compact manifolds. Canad. J. Math. {\bf 52} (5) (2000), 982-998.
\bibitem[N]{N}
T. Napier, Convexity properties of coverings of smooth projective varieties.
Math. Ann. {\bf 286} (1990), 433-479.
\end{thebibliography}
\end{document}